\documentclass[12pt,headings=normal]{scrartcl}

\usepackage[utf8]{inputenc}
\usepackage[ngerman]{babel}
\usepackage{amsmath,amssymb}
\usepackage[thmmarks]{ntheorem}
\usepackage{verbatim}
\usepackage{graphicx}
\usepackage[all]{xy}

\newcommand{\Ass}{\operatorname{Ass}}
\newcommand{\Ann}{\operatorname{Ann}}

\newcommand{\Hom}{\operatorname{Hom}}

\theorembodyfont{\itshape}
\theoremstyle{plain}
\newtheorem{Satz}{Satz}[section]
\newtheorem{Folgerung}[Satz]{Folgerung}
\newtheorem{Lemma}[Satz]{Lemma}

\theorembodyfont{\upshape}
\newtheorem{Bemerkung}[Satz]{Bemerkung}

\theoremheaderfont{\scshape}
\theorembodyfont{\upshape}
\theoremstyle{nonumberplain}
\theoremseparator{}
\theoremsymbol{$\Box$}
\newtheorem{Beweis}{Beweis}

\title{\Large Einfach-teilbare und einfach-torsionsfreie $R$-Moduln}
\author{\large Helmut Zöschinger\\
  \large Mathematisches Institut der Universität München\\
  \large Theresienstr. 39, D-80333 München\\
  \large E-mail: zoeschinger$@$mathematik.uni-muenchen.de
}
\date{}

\newcounter{abccount}
\newenvironment{abc}{%
\begin{list}{(\alph{abccount})}{%
\usecounter{abccount}%
\setlength{\partopsep}{0pt}%
\setlength{\topsep}{1ex}%
\setlength{\itemsep}{0pt}%
}%
}{\end{list}}

\newcounter{myenumcount}

\newcounter{iiicount}
\newenvironment{iii}{%
\begin{list}{(\roman{iiicount})}{%
\usecounter{iiicount}%
\setlength{\labelwidth}{3em}%
\setlength{\partopsep}{0pt}%
\setlength{\topsep}{1ex}%
\setlength{\itemsep}{0pt}%
}%
}{\end{list}}

\newcounter{iiiscount}
\newenvironment{iiis}{%
\begin{list}{(\roman{iiiscount}')}{%
\usecounter{iiiscount}%
\setlength{\labelwidth}{3em}%
\setlength{\partopsep}{0pt}%
\setlength{\topsep}{1ex}%
\setlength{\itemsep}{0pt}%
}%
}{\end{list}}

\begin{document}
\maketitle

\centerline{\textbf{Abstract}}
\begin{abstract}
  \noindent
  Let $(R, \mathfrak{m})$ be a commutative Noetherian local ring with total
  quotient ring $K$. An $R$-module $M$ is called simple divisible, if $M$ is
  divisible $\neq 0$, but every proper submodule $0 \neq U \subsetneqq M$ is
  not divisible. Dually, $M$ is called simple torsion free, if $M$ ist
  torsion free $\neq 0$, but, for every proper submodule $0 \neq U
  \subsetneqq M$, the factor module $M/U$ is not torsion free. Our first
  result is that $M \neq 0$ is simple torsion free iff $M$ is a submodule of
  $\kappa(\mathfrak{p}) = R_{\mathfrak{p}}/\mathfrak{p} R_{\mathfrak{p}}$
  for a maximal element $\mathfrak{p}$ in $\operatorname{Ass}(R)$. The
  structure of simple divisible modules is more complicated and was examined
  primarily by E. Matlis (1973) over 1-dimensional local
  $CM$-rings and by A. Facchini (1989) over any integral domain.
  Our main results are: If the injective hull $E(R/\mathfrak{q})$ is simple
  divisible ($\mathfrak{q} \in \operatorname{Spec}(R)$), then the ring
  $R_{\mathfrak{q}}$ is analytically irreducible and essentially complete.
  Especially for $\mathfrak{q} = \mathfrak{m}$, the simple divisible
  submodules of $E(R/\mathfrak{m})$ correspond exactly to the maximal ideals
  of the ring $\widehat{R} \otimes_R K$, and $E(R/\mathfrak{m})$ itself is
  simple divisible iff $\widehat{R} \otimes_R K$ is a field.
\end{abstract}

\medskip


\bigskip


\section{Prime und koprime $R$-Moduln}

Stets sei $R$ ein kommutativer, noetherscher, lokaler Ring, $\mathfrak{m}$
sein maximales Ideal, $k = R/\mathfrak{m}$ sein Restklassenkörper, $E$ die
injektive Hülle von $k$ und $M^\circ = \Hom_R(M,E)$ das Matlis-Duale eines
$R$-Moduls $M$. Bekanntlich heißt $M$ \emph{prim}, wenn $M \neq 0$ ist und
aus $0\neq U \subset M$ folgt $\Ann_R(U) = \Ann_R(M)$. Das ist äquivalent
damit, dass $M \neq 0$ ist und für jedes $r \in R$ die Multiplikation mit $r
\colon M \to M$ entweder Null oder injektiv ist. Dual heißt ein $R$-Modul
$N$ \emph{koprim} (second, dual prime), wenn $N \neq 0$ ist und aus $V
\subsetneqq N$ folgt $\Ann_R(N/V) =  \Ann_R(N)$, d.\,h. für jedes $r \in R$
die Multiplikation mit $r \colon N \to N$ entweder Null oder surjektiv ist.

\bigskip

\noindent
\textbf{Beispiel 1} Ist $R$ ein Integritätsring, so ist jeder torsionsfreie
$R$-Modul $M \neq 0$ prim und jeder teilbare $R$-Modul $N \neq 0$ koprim.\hfill$\Box$

\medskip

\noindent
\textbf{Beispiel 2} Für jedes Primideal $\mathfrak{p} \in
\operatorname{Spec}(R)$ ist $\kappa(\mathfrak{p}) = R_{\mathfrak{p}} /
\mathfrak{p} R_{\mathfrak{p}}$ als $R$-Modul sowohl prim als auch koprim.

\begin{Beweis}
  Falls $r \in \mathfrak{p}$, ist die Multiplikation mit $r \colon
  \kappa(\mathfrak{p}) \to \kappa(\mathfrak{p})$ die Nullabbildung: $r
  \left[\frac{a}{s}\right] = \left[\frac{ra}{s}\right] = 0$.\\
  Falls $r \notin \mathfrak{p}$, ist sie injektiv:
  $r\left[\frac{a}{s}\right] = 0$ bedeutet $\frac{ra}{s} \in
  \mathfrak{p}R_{\mathfrak{p}}$, $a \in \mathfrak{p}$,
  $\left[\frac{a}{s}\right] = 0$, aber auch surjektiv:
  $\left[\frac{b}{t}\right] = r\left[\frac{b}{rt}\right]$.
\end{Beweis}

\medskip

\noindent
\textbf{Beispiel 3} $M$ prim ${}\iff M^\circ$ koprim. $N$ koprim ${}\iff
N^\circ$ prim.

\begin{Beweis}
  Für jeden $R$-Modul $X$ ist $f \colon X \to X$ genau dann injektiv
  (surjektiv), wenn $f^\circ \colon X^\circ \to X^\circ$ surjektiv
  (injektiv) ist.
\end{Beweis}

\begin{Lemma}\label{1.1}
  (a) Ist $M$ prim, so ist $\Ann_R(M)$ ein Primideal, sagen wir
  $\mathfrak{p}$, und damit gilt $\Ass(M) = \operatorname{Koatt}(M) =
  \{\mathfrak{p}\}$.\\
  (b) Ist $N$ koprim, so ist $\Ann_R(N)$ ein Primideal, sagen wir
  $\mathfrak{p}$, und damit gilt $\operatorname{Koass}(N) =
  \operatorname{Att}(N) = \{p\}$.
\end{Lemma}

\begin{Beweis}
  (a) Klar ist $\Ann_R(M) \neq R$, und aus $r_1 \notin \Ann_R(M) \not\ni
  r_2$ folgt, dass die Multiplikation mit $r_1$ und $r_2$ injektiv ist, also
  auch mit $r_1r_2 \implies r_1r_2 \notin \Ann_R(M)$. Weiter gilt in
  $\varnothing \neq \Ass(M) \subset \operatorname{Koatt}(M)$ für jedes
  $\mathfrak{q} \in \operatorname{Koatt}(M)$, d.\,h. $\mathfrak{q} =
  \Ann_R(M[\mathfrak{q}])$, dass nach Voraussetzung $\Ann_R(M[\mathfrak{q}])
  = \Ann_R(M)$ ist, also $\mathfrak{q} = \mathfrak{p}$.\\
  (b) $N^\circ$ ist nach Beispiel~3 prim, so dass wegen
  $\operatorname{Koass}(N) = \Ass(N^\circ)$, $\operatorname{Att}(N) =
  \operatorname{Koatt}(N^\circ)$ alles mit (a) folgt.
\end{Beweis}

\begin{Satz}\label{1.2}
  Für ein Primideal $\mathfrak{q} \in \operatorname{Spec}(R)$ sind
  äquivalent:
  \begin{iii}
    \item $M = R_{\mathfrak{q}}$ ist als $R$-Modul prim.
    \item $N = E(R/\mathfrak{q})$ ist als $R$-Modul koprim.
    \item $R_{\mathfrak{q}}$ ist ein Integritätsring.
    \item $\mathfrak{p} = \Ann_R(M)$ ist ein Primideal und $M
      \hookrightarrow \kappa(\mathfrak{p})$.
  \end{iii}
  Erfüllt $\mathfrak{q}$ diese vier äquivalenten Bedingungen, ist
  $\mathfrak{p} = \Ann_R(M)$ sogar ein Element von $\Ass(R)$ und
  $\Ann_R(\mathfrak{p}) \not\subset \mathfrak{q}$.
\end{Satz}

\begin{Beweis}
  (i $\leftrightarrow$ ii) Nach dem Beweis von (\cite{010} Lemma~2.3) ist
  $N$ ein direkter Summand von $M^\circ$ und $M$ ein Untermodul von
  $N^\circ$ (auch bei $\mathfrak{q} = \mathfrak{m})$, so dass Beispiel~3 die
  Behauptung liefert.\\
  (i $\leftrightarrow$ iii) Für jedes $r \in R$ gilt: Genau dann ist $rM =
  0$, wenn $\frac{r}{1} = 0$ ist im Ring $R_{\mathfrak{q}}$; genau dann ist
  $M[r] = 0$, wenn $\frac{r}{1}$ ein NNT im Ring $R_{\mathfrak{q}}$ ist.
  Damit folgt die Behauptung.\\
  (iii $\to$ iv) \emph{1. Schritt} $M = R_{\mathfrak{q}}$ ist als $R$-Modul
  uniform: Die Elemente von $S = R \setminus \mathfrak{q}$ operieren auf $M$
  bijektiv, so dass für jeden $R$-Untermodul $0 \neq U \subset M$ auch $0
  \neq U_S \subset M_S$ ist, mit $0 \neq V \subset M$ ebenso $0 \neq V_S
  \subset M_S$, also (weil $M_S$ als $R_S$-Modul uniform ist) $U_S \cap V_S
  \neq 0$, $U \cap V \neq 0$. \emph{2. Schritt} $\mathfrak{p} = \Ann_R(M)$
  ist nach (\ref{1.1}) das einzige zu $M$ assoziierte Primideal, also $M
  \subset E(R/\mathfrak{p})[\mathfrak{p}] \cong \kappa(\mathfrak{p})$.\\
  (iv $\to$ i) $\kappa(\mathfrak{p})$ ist nach Beispiel~2 prim, also auch
  der Untermodul $M$.\\
  Der Zusatz gilt, weil $M$ flach, also $\Ass(M) \subset \Ass(R)$ ist, und
  aus $\mathfrak{p}M = 0$, $\mathfrak{p}R_{q} = 0$ folgt auch
  $\Ann_R(\mathfrak{p}) \not\subset \mathfrak{q}$. 
\end{Beweis}

\begin{Bemerkung}\label{1.3}
  Aus der Bedingung (iv) folgt, weil $\kappa(\mathfrak{p})^\circ \cong
  \kappa(\mathfrak{p})^{(I)}$ und $N$ direkter Summand von $M^\circ$ ist,
  dass es einen Epimorphismus
  \begin{equation*}
    \kappa(\mathfrak{p})^{(I)} \twoheadrightarrow N
  \end{equation*}
  gibt. Man kann aber \emph{nicht} $|I|=1$ erwarten: Ist $R$ ein
  1-dimensionaler Integritätsring mit Quotientenkörper $K$ und $\mathfrak{q}
  = \mathfrak{m}$ (also $\mathfrak{p}=0$), so gibt es nach Matlis
  (\cite{004} Theorem~15.7) genau dann einen Epimorphismus $K
  \twoheadrightarrow E$, wenn $\widehat{R} \otimes_R K$ als
  $\widehat{R}$-Modul injektiv ist (d.\,h. $R$ ein kanonisches Ideal
  besitzt).\hfill$\Box$
\end{Bemerkung}

\begin{Bemerkung}\label{1.4}
  Ist $R_{\mathfrak{q}}$ in der Bedingung (iii) sogar ein Körper (d.\,h.
  $\mathfrak{q}R_{\mathfrak{q}} = 0$), folgt $M \cong N \cong
  \kappa(\mathfrak{q})$. Man kann aber \emph{nicht} erwarten, dass
  $\kappa(\mathfrak{q})$ einfach-torsionsfrei oder einfach-teilbar ist:
  Wählt man $R$ und $\mathfrak{q}$ so, dass $0 \neq \mathfrak{q}
\subsetneqq \mathfrak{m}$ und $\operatorname{So}(R) =\mathfrak{q}$ ist, wird
(wegen $\operatorname{So}(R) \neq 0$) jeder $R$-Modul torsionsfrei und
teilbar, aber nur $k \cong \kappa(\mathfrak{m})$ einfach-torsionsfrei
(einfach-teilbar).\hfill$\Box$
\end{Bemerkung}

\section{Die Struktur der einfach-torsionsfreien $R$-Moduln}

\begin{Lemma}\label{2.1}
  Ist $M$ einfach-torsionsfrei, so gilt:
  \begin{abc}
  \item Für jeden Untermodul $0 \neq U \subset M$ ist $M/U$ torsion und
    $U$ wieder einfach-torsionsfrei.
  \item $M$ ist uniform.
  \item $\Ann_R(M)$ ist ein maximales Element in $\Ass(R)$.
  \end{abc}
\end{Lemma}

\begin{Beweis}
  (a) Mit $U_*/U = \operatorname{T}(M/U)$ ist $M/U_*$ torsionsfrei, also
  nach Voraussetzung $U_* = M$, d.\,h. $M/U$ torsion. Insbesondere folgt aus
  $0 \neq U_1 \subsetneqq U$, dass auch $U/U_1$ torsion ist, also nicht
  torsionsfrei wie behauptet.\\
  (b) Aus $0 \neq U_1 \subset M$, $0 \neq U_2 \subset M$ folgt $M/U_1$,
  $M/U_2$ torsion, also auch $M/U_1 \cap U_2$, und weil $M$ torsionsfrei
  ist, $U_1 \cap U_2 \neq 0$.\\
  (c) Nach Lemma~\ref{1.1} ist $\mathfrak{p} = \Ann_R(M)$ ein Primideal, das
  (weil $M$ torsionsfrei ist) keinen NNT enthält, d.\,h. nicht regulär ist.
  Für jedes echt größere Primideal $\mathfrak{p} \subsetneqq \mathfrak{p}_1$
  ist aber $R/\mathfrak{p}_1$ torsion (weil $R/\mathfrak{p} \hookrightarrow
  M$ nach (a) wieder einfach-torsionsfrei ist), also $\mathfrak{p}_1$ regulär.
\end{Beweis}

\noindent
\textbf{Beispiel 1} Sei $R$ ein Integritätsring. Genau dann ist $M$
einfach-torsionsfrei, wenn $M$ torsionsfrei und uniform ist.

\begin{Beweis}
  "`$\Rightarrow$"' ist (b), und bei "`$\Leftarrow$"' folgt mit einer
  injektiven Hülle $M \subset X$, dass auch $X$ torsionsfrei und uniform
  ist, also isomorph zum Quotientenkörper $K$. Natürlich ist $K$ als
  $R$-Modul einfach-torsionsfrei, also auch $M$.
\end{Beweis}

\noindent
\textbf{Beispiel 2} Sei $\mathfrak{p} \in \operatorname{Spec}(R)$. Genau
dann ist $R/\mathfrak{p}$ als $R$-Modul einfach-torsionsfrei, wenn
$\mathfrak{p}$ ein maximales Element in $\Ass(R)$ ist.

\begin{Beweis}
  "`$\Rightarrow$"' ist (c), und bei "`$\Leftarrow$"' gilt für jedes Ideal
  $\mathfrak{p} \subsetneqq \mathfrak{a} \subsetneqq R$, dass in
  $\sqrt{\mathfrak{a}} = \mathfrak{q}_1 \cap \ldots \cap \mathfrak{q}_n$
  alle $\mathfrak{q}_i \not\subset \bigcup \Ass(R)$, also regulär sind,
  so dass auch $\mathfrak{a}$ regulär ist, d.\,h. $R/\mathfrak{a}$ torsion.
\end{Beweis}

\begin{Satz}\label{2.2}
  Für einen $R$-Modul $M \neq 0$ sind äquivalent:
  \begin{iii}
    \item $M$ ist einfach-torsionsfrei.
    \item $M \hookrightarrow \kappa(\mathfrak{p})$ für ein maximales Element
      $\mathfrak{p}$ in $\Ass(R)$.
  \end{iii}
\end{Satz}

\begin{Beweis}
  (i $\to$ ii) Nach (\ref{2.1}) ist $\mathfrak{p} = \Ann_R(M)$ ein maximales
  Element in $\Ass(R)$, außerdem $M$ uniform mit dem nach (\ref{1.1})
  einzigen assoziierten Primideal $\mathfrak{p}$. Es folgt $M \hookrightarrow
  E(R/\mathfrak{p})[\mathfrak{p}] \cong \kappa(\mathfrak{p})$.\\
  (ii $\to$ i) Wieder nach (\ref{2.1}) vererbt sich die Eigenschaft
  "`einfach-torsionsfrei"' auf Untermoduln $\neq 0$, so dass wir sie nur noch
  für $\kappa(\mathfrak{p})$ zeigen müssen. \emph{1.~Schritt} Wegen der
  Maximalität von $\mathfrak{p}$ ist $\operatorname{Kok}(R/\mathfrak{p}
  \overset{\varepsilon}{\hookrightarrow} \kappa(\mathfrak{p}))$ $R$-torsion,
  denn für jedes $[\frac{a}{s}] \in \kappa(\mathfrak{p})$ ist $\mathfrak{p}
  + (s)$ regulär, also $r - bs \in \mathfrak{p}$ für ein $b \in R$ und einen
  NNT $r \in R$, so dass $r [\frac{a}{s}] = [\frac{ba}{1}] \in
  \operatorname{Bi} \varepsilon$ folgt. \emph{2.~Schritt} Weil
  $\mathfrak{p}$ nicht regulär ist, ist $\kappa(\mathfrak{p})$ als $R$-Modul
  torsionsfrei (vgl. Beispiel~2 in (1)), und ist jetzt $U \subset
  \kappa(\mathfrak{p})$ ein $R$-Untermodul $\neq 0$, folgt aus $U \cap
  \operatorname{Bi} \varepsilon \neq 0$ nach Beispiel~2 $\operatorname{Bi}
  \varepsilon / U \cap \operatorname{Bi} \varepsilon$ $R$-torsion, also mit
  dem 1.~Schritt sogar $\kappa(\mathfrak{p})/U$ $R$-torsion wie gewünscht. 
\end{Beweis}

\begin{Folgerung}\label{2.3}
  Für einen $R$-Modul sind äquivalent:
  \begin{iii}
    \item $M$ ist einfach-torsionsfrei und teilbar.
    \item $M$ ist torsionsfrei und einfach-teilbar.
    \item $M \cong \kappa(\mathfrak{p})$ für ein maximales Element
      $\mathfrak{p}$ in $\Ass(R)$.
  \end{iii}
\end{Folgerung}

\begin{Beweis}
  (i $\to$ iii) Nach dem Satz ist $M \cong M' \subset \kappa(\mathfrak{p})$
  mit $\kappa(\mathfrak{p})$ einfach-torsionsfrei. Weil $M'$ teilbar, also
  nach dem Schlangenlemma auch $\kappa(\mathfrak{p})/M'$ torsionsfrei ist,
  folgt mit (\ref{2.1}, a) bereits $M' = \kappa(\mathfrak{p})$.\\
  (iii $\to$ ii) Es ist nur noch zu zeigen, dass $\kappa(\mathfrak{p})$
  einfach-teilbar ist: Klar ist es teilbar (vgl. Beispiel~2 in (1)), und aus
  $0 \neq V \subset \kappa(\mathfrak{p})$, $V$ teilbar, folgt nach dem
  Schlangenlemma $\kappa(\mathfrak{p})/V$ torsionsfrei, also wieder nach
  (\ref{2.1}, a) $V = \kappa(\mathfrak{p})$.\\
  (ii $\to$ i) Es ist nur noch zu zeigen, dass $M$ einfach-torsionsfrei ist:
  Aus $0 \neq U \subset M$ und $M/U$ torsionsfrei folgt wieder nach dem
  Schlangenlemma, dass $U$ teilbar, also $U=M$ ist.
\end{Beweis}

\section{Einfach-teilbare $R$-Moduln}
Eine dem Satz~(\ref{2.2}) entsprechende Strukturaussage für einfach-teilbare
$R$-Moduln gibt es nicht. Die in (\ref{2.1}) gezeigten Eigenschaften gelten
im dualen Fall nicht mehr, und selbst bei $N = E(R/\mathfrak{q})$ können wir
nur für spezielle Primideale $\mathfrak{q}$ entscheiden, wann $N$
einfach-teilbar ist.

\begin{Lemma}\label{3.1}
  Ist $N$ einfach-teilbar, so gilt:
  \begin{abc}
  \item $N$ ist torsionsfrei oder torsion.
  \item Ist $N$ komplementiert, so ist $N$ kouniform.
  \item Hat der Untermodul $V \subsetneqq N$ in jedem Zwischenmodul ein
    Komplement, so ist auch $N/V$ einfach-teilbar.
  \item Ist der Untermodul $V \subsetneqq N$ artinsch, so ist $V$
    beschränkt und es gibt einen Epimorphismus $N/V \twoheadrightarrow N$.
  \end{abc}
\end{Lemma}

\begin{Beweis}
  (a) In jedem teilbaren $R$-Modul $X$ ist $\operatorname{T}(X)$ wieder
  teilbar, in unserem Fall also $\operatorname{T}(N) = 0$ oder
  $\operatorname{T}(N) = N$.\\
  (b) Jeder echte Untermodul $V \subsetneqq N$ hat nach Voraussetzung ein
  Komplement $W$ in $N$, und weil $W$ teilbar $\neq 0$ ist, folgt $W=N$, $V$
  klein in $N$.\\
  (c) Klar ist $N/V$ teilbar $\neq 0$. Gäbe es einen teilbaren Untermodul $0
  \neq W/V \subsetneqq N/V$, folgte mit einem Komplement $W'$ von $V$ in
  $W$, dass auch $W'$ teilbar wäre und $0 \neq W' \subsetneqq N$ entgegen
  der Voraussetzung.\\
  (d) Die Menge $\{rV \,|\, r \in R \text{ ein NNT}\}$ hat ein minimales
  Element $r_0 V$, und weil $r_0V$ teilbar ist, gilt schon $r_0V=0$. Weiter
  folgt aus $V \subset N[r_0]$ der Epimorphismus $N/V \twoheadrightarrow
  N/N[r_0] \cong r_0 N = N$. 
\end{Beweis}

\begin{Lemma}\label{3.2}
  \begin{abc}
  \item $N^\circ$ einfach-torsionsfrei $\implies$ $N$ einfach-teilbar.
  \item $N^\circ$ einfach-teilbar $\implies$ $N$ einfach-torsionsfrei.
  \item Ist $N$ reflexiv, gilt beide Male die Umkehrung.
  \end{abc}
\end{Lemma}

\begin{Beweis}
  (a) Klar ist $N$ teilbar $\neq 0$. Ist aber $0 \neq V \subsetneqq N$,
  folgt aus $0 \neq \Ann_{N^\circ}(V) \subsetneqq N^\circ$, dass
  $N^\circ/\Ann_{N^\circ}(V) \cong V^\circ$ nicht torsionsfrei ist, also $V$
  nicht teilbar.\\
  (b) Entsprechend.\\
  (c) Folgt mit $N \cong N^{\circ\circ}$ aus (a) bzw. (b). 
\end{Beweis}

\begin{Folgerung}\label{3.3}
  Ist $N$ einfach-teilbar und \emph{reflexiv}, so gilt:
  \begin{abc}
  \item Für jeden Untermodul $V \subsetneqq N$ ist $V$ kotorsion und $N/V$
    wieder einfach-teilbar.
  \item $N$ ist kouniform.
  \item Es gibt einen Epimorphismus $\kappa(\mathfrak{p}) \twoheadrightarrow
    N$ mit einem maximalen Element $\mathfrak{p}$ in $\Ass(R)$.
  \end{abc}
\end{Folgerung}

\begin{Beweis}
  (a) Aus $(N/V)^\circ \hookrightarrow N^\circ$ folgt mit (\ref{2.1}, a),
  dass $(N/V)^\circ$ einfach-torsionsfrei ist und $V^\circ$ torsion, also
  $N/V$ einfach-teilbar ist und $V$ kotorsion (d.\,h. jeder artinsche
  Faktormodul von $V$ beschränkt).\\
  (b) Nach (\ref{2.1}, b) ist $N^\circ$ uniform, also $N$ kouniform.\\
  (c) Nach (\ref{2.2}) gibt es ein maximales Element $\mathfrak{p}$ in
  $\Ass(R)$ mit $N^\circ \hookrightarrow \kappa(\mathfrak{p})$, so dass aus
  $\kappa(\mathfrak{p})^\circ \cong \kappa(\mathfrak{p})^{(I)}$ und
  $\kappa(\mathfrak{p})^\circ \twoheadrightarrow N$ die Behauptung folgt.
\end{Beweis}

Weil in einem reflexiven Modul jeder Untermodul komplementiert ist, folgen
die Punkte (a) und (b) auch aus~(\ref{3.1}).

\bigskip

\noindent
\textbf{Beispiel 1} Sei $R$ ein Integritätsring mit Quotientenkörper $K$ und
einem Primideal $0 \neq \mathfrak{q} \subsetneqq \mathfrak{m}$. Dann ist $X
= R_{\mathfrak{q}}$ eine nichttrivialer Untermodul von $K$ und $X$ nach
(\cite{006} Lemma~4.3) \emph{nicht} klein in $K$, insbesondere nicht
kotorsion. Aus $Y \subsetneqq K$, $Y+X=K$ folgt auch, dass $K$ nicht
kouniform ist und $K/Y\cap X$ nicht einfach-teilbar.\hfill$\Box$

\bigskip

Wir wollen untersuchen, wann $N = E(R/\mathfrak{q})$ als $R$-Modul
einfach-teilbar ist. Der Spezialfall $\mathfrak{q} = \mathfrak{m}$ wurde
bereits in (\cite{009} Lemma~4.3) behandelt:

\medskip
\noindent
\textbf{Beispiel 2} Genau dann ist $E$ als $R$-Modul einfach-teilbar, wenn
die Vervollständigung $\widehat{R}$ ein Integritätsring ist und wenn aus $Q
\in \operatorname{Spec}(\widehat{R})$, $Q \cap R = 0$ stets folgt $Q =
0$.\hfill$\Box$

Die erste Bedingung bedeutet, dass $R$ analytisch irreduzibel ist, die
zweite, dass $R$ im Sinne von (\cite{007} p. 3402) im wesentlichen
vollständig ist.

\bigskip

\noindent
\textbf{Beispiel 3} Ist $\mathfrak{q} \in \operatorname{Spec}(R)$
\emph{nicht} regulär, so sind äquivalent:
\begin{iiis}
  \item $M = R_{\mathfrak{q}}$ ist einfach-torsionsfrei.
  \item $N = E(R/\mathfrak{q})$ ist einfach-teilbar.
  \item $R_{\mathfrak{q}}$ ist ein Körper und $\mathfrak{q}$ ein maximales
    Element in $\Ass(R)$.
\end{iiis}

\begin{Beweis}
  (i' $\to$ iii') Nach (\ref{2.1}) ist $\mathfrak{p} = \Ann_R(M)$ ein
  maximales Element in $\Ass(R)$, so dass aus $\mathfrak{p} \subset
  \mathfrak{q} \subset \mathfrak{q}_1 \in \Ass(R)$ folgt $\mathfrak{p} =
  \mathfrak{q}$, also die zweite Aussage bewiesen ist. Wegen
  $\mathfrak{q}R_{\mathfrak{q}} = 0$ gilt auch die erste.\\
  (ii' $\to$ iii') Klar ist $R/\mathfrak{q}$ torsionsfrei, also auch $N$,
  so dass $N$ nach (\ref{2.3}) sogar einfach-torsionsfrei ist und
  $\mathfrak{p} = \Ann_R(M)$ ein maximales Element in $\Ass(R)$. Wie eben
  folgt, dass $\mathfrak{p} = \mathfrak{q}$ und $R_{\mathfrak{q}}$ ein
  Körper ist.\\
  Gilt aber (iii'), ist $M \cong N \cong \kappa(\mathfrak{q})$ wieder nach
  (\ref{2.3}) sowohl einfach-torsionsfrei als auch einfach-teilbar.
\end{Beweis}

Hat z.\,B. $R$ keine nilpotenten Elemente, gelten die drei äquivalenten
Aussagen für jedes $\mathfrak{q} \in \Ass(R)$.

\begin{Lemma}\label{3.4}
  Sei $N$ ein $R$-Modul und $\mathfrak{q} \in \operatorname{Spec}(R)$,
  so dass alle $s \in R \setminus \mathfrak{q}$ auf $N$ bijektiv operieren.
  Ist dann noch $R \to R_{\mathfrak{q}}$ injektiv, so sind äquivalent:
  \begin{iii}
    \item $N$ ist als $R$-Modul einfach-teilbar.
    \item $N$ ist als $R_{\mathfrak{q}}$-Modul einfach-teilbar.
  \end{iii}
\end{Lemma}

\begin{Beweis}
  (i $\to$ ii) $N$ ist als $R_{\mathfrak{q}}$-Modul teilbar, denn für jeden
  NNT $\frac{r}{s} \in R_{\mathfrak{q}}$ ist nach der zweiten Voraussetzung
  auch $r \in R$ ein NNT, also $rN=N$, $\frac{r}{s} \cdot N = N$. Gäbe es
  einen $R_{\mathfrak{q}}$-Untermodul $0 \neq X \subsetneqq N$, der als
  $R_{\mathfrak{q}}$-Modul teilbar ist, folgte für alle NNT $r \in R$, dass
  $\frac{r}{1} \cdot X = X$, also $rX=X$ ist, d.\,h. $X$ als $R$-Modul
  teilbar entgegen der Annahme.\\
  (ii $\to$ i) $N$ ist als $R$-Modul teilbar, denn für jeden NNT $r \in R$
  ist $\frac{r}{1}$ ein NNT in $R_{\mathfrak{q}}$, so dass aus $\frac{r}{1}
  \cdot N = N$ folgt $rN=N$. Gäbe es einen $R$-Untermodul $0 \neq V
  \subsetneqq N$, der als $R$-Modul teilbar ist, wäre $V$ auch ein
  $R_{\mathfrak{q}}$-Untermodul von $N$ ($\frac{r}{s} \in R_{\mathfrak{q}}
  \implies s$ ein NNT in $R$, $sV=V$, $\frac{r}{s} \cdot V \subset V$) und
  als $R_{\mathfrak{q}}$-Modul teilbar ($\frac{r}{s} \in R_{\mathfrak{q}}$
  ein NNT $\implies r \in R$ ein NNT, $sV \subset rV$, $V \subset
  \frac{r}{s} \cdot V$) entgegen der Annahme.
\end{Beweis}

\begin{Satz}\label{3.5}
  Für ein Primideal $\mathfrak{q} \in \operatorname{Spec}(R)$ betrachten wir
  folgende Eigenschaften:
  \begin{iii}
    \item $E(R/\mathfrak{q})$ ist als $R$-Modul einfach-teilbar.
    \item Der Ring $R_{\mathfrak{q}}$ ist analytisch irreduzibel und im
      wesentlichen vollständig.
  \end{iii}
  Dann gilt (i $\to$ ii), und falls $R \to R_{\mathfrak{q}}$ injektiv ist,
  auch (ii $\to$ i).  
\end{Satz}

\begin{Beweis}
  Natürlich operieren alle $s \in R \setminus \mathfrak{q}$ auf $N =
  E(R/\mathfrak{q})$ bijektiv, und als $R_{\mathfrak{q}}$-Modul ist $N$ die
  injektive Hülle von $\kappa(\mathfrak{q})$.\\
  Bei (i $\to$ ii) ist ohne weitere Voraussetzung $N$ als
  $R_{\mathfrak{q}}$-Modul teilbar, ja nach dem Beweis in (\ref{3.4}) sogar
  einfach-teilbar, so dass Beispiel~2 die Behauptung liefert.\\
  Bei (ii $\to$ i) ist wieder nach dem Beispiel~2 $N$ als
  $R_{\mathfrak{q}}$-Modul einfach-teilbar, mit der Zusatzbedingung und
  (\ref{3.4}) also (i) erfüllt.
\end{Beweis}

In (\ref{1.4}) ist $R_{\mathfrak{q}}$ ein Körper (also (ii) erfüllt), aber
$E(R/\mathfrak{q})$ nicht einfach-teilbar.

\begin{Folgerung}\label{3.6}
  Ist $E(R/\mathfrak{q})$ einfach-teilbar und $R$ abzählbar, folgt
  $\operatorname{h}(\mathfrak{q}) \leq 1$.
\end{Folgerung}

\begin{Beweis}
  Nach dem Satz genügt es zu zeigen: Ist $R$ analytisch irreduzibel, im
  wesentlichen vollständig und abzählbar, folgt $\dim(R) \leq 1$. Wäre
  $\dim(R) \geq 2$, existierte nach Heinzer--Rotthaus--Sally (\cite{003}
  Proposition~4.10) ein $Q \in \operatorname{Spec}(\widehat{R})$ mit
  $\dim(\widehat{R}/Q) = 1$ und $Q \cap R = 0$ im Widerspruch zur zweiten
  Bedingung.
\end{Beweis}

\begin{Bemerkung}\label{3.7}
  Sei $K$ der totale Quotientenring von $R$ und $\mathfrak{N}$ das
  Nilradikal von $\widehat{R}$. Betrachten wir die folgenden vier Bedingungen
  \begin{abc}
  \item Der Ring $\widehat{R} \otimes_R K$ ist halbeinfach,
  \item Es ist $\mathfrak{N} = 0$ und $\widehat{R} \otimes_R K$ als
    $\widehat{R}$-Modul injektiv (schwach-injektiv, teilbar),
  \item $R$ ist im wesentlichen vollständig, d.\,h. aus $A$ groß in
    $\widehat{R}$ folgt $A \cap R$ groß in $R$,
  \item Aus $P \in \operatorname{Spec}(\widehat{R})$, $P$ groß in
    $\widehat{R}$ folgt $P \cap R$ groß in $R$,
  \end{abc}
  so gilt stets (a $\iff$ b $\implies$ c $\implies$ d), und falls $R$ keine
  nilpotenten Elemente hat, auch (d $\implies$ b).\hfill$\Box$
\end{Bemerkung}

Die Bedingung (a) führt zu einer weiteren Beschreibung von Beispiel~2: Genau
dann ist $E$ als $R$-Modul einfach-teilbar, wenn $\widehat{R} \otimes_R K$
ein Körper ist.

Die Bedingung (b) ist nach (\cite{008} Lemma~2.12) genau dann erfüllt, wenn
jeder Faktormodul von ${}_RK$ schwach-injektiv ist, d.\,h. $\overline{Z}(K)
= K$.

\section{Die Struktur der einfach-teilbaren artinschen $R$-Moduln}

Jeder halbartinsche $R$-Modul $N$ lässt sich in natürlicher Weise zu einem
$\widehat{R}$-Modul machen und wir wollen zunächst prüfen, wann ein
Untermodul $0 \neq N \subset E$ als $R$-Modul oder als $\widehat{R}$-Modul
einfach-teilbar (teilbar, koprim) ist.

\begin{Lemma}\label{4.1}
  Sei $N$ als $R$-Modul einfach-teilbar und artinsch. Dann gilt:
  \begin{abc}
  \item $N$ ist als $\widehat{R}$-Modul koprim.
  \item $N$ ist äquivalent zu genau einem einfach-teilbaren Untermodul $N_1
    \subset E$.
  \end{abc}
\end{Lemma}

\begin{Beweis}
  (a) Ist $u \in\widehat{R}$ und $u \cdot N \neq 0$, so ist der
  $R$-Homomorphismus $N \to N$, $y \mapsto u \cdot y$ ungleich Null, also
  nach Voraussetzung surjektiv, d.\,h. es folgt $u \cdot N = N$.\\
  (b) Existenz: Es gibt einen Untermodul $V \subsetneqq N$ mit $N/V
  \hookrightarrow E$, und nach (\ref{3.1}, d) folgt $N/V \sim N$.\\
  Eindeutigkeit: Zwei äquivalente Untermoduln $X,Y$ von $E$ sind bereits
  identisch, denn allein aus $\alpha\colon X \twoheadrightarrow Y$ folgt,
  dass $\alpha$ auch $\widehat{R}$-linear ist, also $\Ann_{\widehat{R}}(X)
  \subset \Ann_{\widehat{R}}(Y)$, $\Ann_E \Ann_{\widehat{R}}(X) \supset
  \Ann_E \Ann_{\widehat{R}}(Y)$, d.\,h. $X \supset Y$.
\end{Beweis}

\begin{Lemma}\label{4.2}
  Sei $0 \neq N \subset E$ und $B = \Ann_{\widehat{R}}(N)$. Dann sind äquivalent:
  \begin{iii}
    \item $N$ ist als $R$-Modul \emph{koprim}.
    \item Für jeden Primdivisor $P$ von $B$ in $\widehat{R}$ gilt $P \cap R
      = B \cap R$.
  \end{iii}
\end{Lemma}

\begin{Beweis}
  Nach (\cite{009} Lemma~4.1) ist $\operatorname{Koass}_R(N) = \{P \cap R
  \,|\, P \text{ ist ein Primdivisor von } B \text{ in } \widehat{R}\}$. Bei
  (i $\to$ ii) ist deshalb $(P \cap R)N = 0$, d.\,h. $P \cap R \subset B
  \cap R$ für diese $P$, und bei (ii $\to$ i) gilt für $rN \neq 0$, d.\,h.
  $r \notin B \cap R$, dass nach Voraussetzung $r$ in keinem $\mathfrak{p}
  \in \operatorname{Koass}(N)$ liegt, also $rN = N$ ist.
\end{Beweis}

Mit denselben Bezeichnungen ist also $N$ als $\widehat{R}$-Modul genau dann
koprim, wenn $B \in \operatorname{Spec}(\widehat{R})$ ist.

\begin{Lemma}\label{4.3}
  Sei $0 \neq N \subset E$ und $B = \Ann_{\widehat{R}}(N)$. Dann sind
  äquivalent:
  \begin{iii}
    \item $N$ ist als $R$-Modul \emph{teilbar}.
    \item Für jeden Primdivisor $P$ von $B$ in $\widehat{R}$ ist $P \cap R$
      nicht regulär.
  \end{iii}
\end{Lemma}

\begin{Beweis}
  Die zweite Bedingung ist äquivalent damit, dass jeder NNT $r \in R$ in
  keinem der genannten $P$, also in keinem $\mathfrak{p} \in
  \operatorname{Koass}_R(N)$ liegt, d.\,h. $rN=N$ ist.
\end{Beweis}

Mit denselben Bezeichnungen ist also $N$ als $\widehat{R}$-Modul genau dann
teilbar, wenn jeder Primdivisor von $B$ in $\widehat{R}$ nicht regulär ist.

\begin{Satz}\label{4.4}
  Sei $0 \neq N \subset E$ und $B = \Ann_{\widehat{R}}(N)$. Dann sind
  äquivalent:
  \begin{iii}
    \item $N$ ist als $R$-Modul \emph{einfach-teilbar}.
    \item $B$ ist ein maximales Elemente in der Menge $\{P \in
      \operatorname{Spec}(\widehat{R}) \,|\,  P\cap R \text{ ist nicht
        regulär}\}$.
  \end{iii}
\end{Satz}

\begin{Beweis}
  (i $\to$ ii) Nach (\ref{4.1}, a) ist $N$ als $\widehat{R}$-Modul
  wenigstens koprim, also $B \in \operatorname{Spec}(\widehat{R})$ nach
  (\ref{1.1}, b) und natürlich $B \cap R$ nicht regulär. Ist aber $B \subset
  P \in \{\}$, wird $V := \Ann_E(P)$ nach (\ref{4.3}) ein teilbarer
  $R$-Modul mit $0 \neq V \subset N$, so dass folgt $V=N$, $B=P$.\\
  (ii $\to$ i) Wieder nach (\ref{4.3}) ist $N = \Ann_E(B)$ als $R$-Modul
  teilbar. Ist aber in $0 \neq V \subset N$ auch $V$ teilbar, folgt mit $B'
  := \Ann_{\widehat{R}}(V)$ und irgendeinem Primdivisor $P_0$ von $B'$ in
  $\widehat{R}$, dass $P_0 \cap R$ nicht regulär ist, also nach
  Voraussetzung $B = P_0$, $B' = P_0$, $V = \Ann_E(B') = N$.
\end{Beweis}

Mit denselben Bezeichnungen ist also $N$ als $\widehat{R}$-Modul genau dann
einfach-teilbar, wenn $B$ ein maximales Element in $\Ass(\widehat{R})$ ist
(was auch mit Beispiel~2 in (2) folgt).

\bigskip

\noindent
\textbf{Beispiel 1} Sei $R$ analytisch irreduzibel, $\dim(R) \geq 2$ und $Q
\in \operatorname{Spec}(\widehat{R})$ mit $\dim(\widehat{R}/Q) = 1$, $Q \cap
R = 0$ (siehe \cite{009} p.~1990). Dann ist $E$ als $\widehat{R}$-Modul
einfach-teilbar, nicht aber als $R$-Modul. Umgekehrt ist $N := \Ann_E(Q)$
als $R$-Modul einfach-teilbar, aber als $\widehat{R}$-Modul nicht einmal
teilbar.\hfill$\Box$

\bigskip

Die Bedingung (ii) im Satz kann man auch durch die maximalen Ideale des
Ringes $\widehat{R} \otimes_R K$ ausdrücken. Sei dazu im Folgenden stets
$\varphi\colon R \to \widehat{R}$ die kanonische Einbettung, $S$ die Menge
\emph{aller} NNT von $R$ und $T = \varphi(S)$. Dann ist $\widehat{R}
\otimes_R K$ ringisomorph zu $\widehat{R}_T$ (via $u \otimes \frac{1}{s}
\mapsto u/\varphi(s)$), und für jedes Ideal $A \subset \widehat{R}$ gilt
\begin{equation*}
  A \cap R \text{ nicht regulär} \iff A \cap T = \varnothing.
\end{equation*}
Jedem maximalen Element in der Menge $\{P \in
\operatorname{Spec}(\widehat{R}) \,|\, P\cap R \text{ ist nicht regulär}\}$
entspricht also ein maximales Ideal im Ring $\widehat{R}_T$, so dass man mit
(\ref{4.4}) erhält:

\begin{Folgerung}\label{4.5}
  Die einfach-teilbaren Untermoduln von $E$ entsprechen genau den maximalen
  Idealen des Ringes $\widehat{R} \otimes_R K$.\hfill$\Box$
\end{Folgerung}

\noindent
\textbf{Beispiel 2} Der Ring $\widehat{R} \otimes_R K$ kann unendlich viele
maximale Ideale haben. Ist z.\,B. $R$ ein \emph{abzählbarer} Integritätsring
mit $\dim(R) \geq 2$, so hat die generische formale Faser von $R$ nach
(\cite{003} Proposition~4.10) überabzählbar viele maximale Elemente $P$, und
für jedes von ihnen gilt $\dim(\widehat{R}/P) = 1$, wegen $P \cap T =
\varnothing$, also $P \cdot \widehat{R}_T \in
\operatorname{Max}(\widehat{R}_T)$.\hfill$\Box$

\medskip

Selbst für einen einfach-teilbaren $R$-Modul $N \subset E$ muss $\Ann_R(N)$
\emph{kein} maximales Element in $\Ass(R)$ sein (vgl.~2.1,~c), es kann sogar
$\Ann_R(N) \notin \Ass(R)$ sein:

\medskip

\noindent
\textbf{Beispiel 3} Sei $R$ abzählbar und sockelfrei, $\mathfrak{p} \in
\operatorname{Spec}(R)$ nicht regulär, aber \emph{kein} Element von
$\Ass(R)$. Dann gibt es einen einfach-teilbaren $R$-Modul $N \subset E$ mit
$\Ann_R(N) = \mathfrak{p}$.

\medskip

Zum Beweis wähle man ein $\mathfrak{p} \subsetneqq \mathfrak{p}_1 \in
\Ass(R)$. Wegen $\mathfrak{m} \notin \Ass(R)$ folgt $\dim(R/\mathfrak{p})
\geq 2$, so dass nach (\cite{003} Proposition~4.10) ein
$\mathfrak{p}\widehat{R} \subset P \in \operatorname{Spec}(\widehat{R})$
existiert mit $\dim(\widehat{R}/P) = 1$ und $P\cap R=\mathfrak{p}$. Dann
leistet $N = \Ann_E(P)$ das Gewünschte, denn klar ist $\Ann_R(N) =
\mathfrak{p}$, nach (\ref{4.4}) aber $N$ auch einfach-teilbar.\hfill$\Box$

\bigskip

Zum Abschluss wollen wir noch eine Klasse von einfach-teilbaren, artinschen
$R$-Moduln $N$ angeben, bei denen $\Ann_R(N)$ doch ein maximales Element in
$\Ass(R)$ ist. Das ist klar, wenn $N$ reflexiv ist (siehe~\ref{3.2}),
speziell wenn $R$ vollständig ist. Es genügt sogar, dass $R \subset
\widehat{R}$ die Bedingung going up erfüllt, und eine Abschwächung davon
wurde von Cuong--Nhan in \cite{001} eingeführt: Ein artinscher $R$-Modul $M$
erfüllt die Bedingung (*), wenn $\operatorname{Koatt}(M) = V(\Ann_R(M))$
ist.

\begin{Lemma}\label{4.6}
  Sei $N$ einfach-teilbar und artinsch mit der Eigenschaft (*). Dann ist
  $\Ann_R(N)$ ein maximales Element in $\Ass(R)$.
\end{Lemma}

\begin{Beweis}
  Sei im \emph{1. Schritt} $N \subset E$. Dann ist $\mathfrak{p} =
  \Ann_R(N)$ ein nicht-reguläres Primideal, und es bleibt zu zeigen, dass
  jedes echt größere Primideal $\mathfrak{p} \subsetneqq \mathfrak{p}_1$
  regulär ist: Mit $P = \Ann_{\widehat{R}}(N)$ bedeutet (*), dass
  $R/\mathfrak{p} \hookrightarrow \widehat{R}/P$ die Bedingung (LO) erfüllt
  (\cite{009} Folgerung~3.2), also ein $P \subsetneqq P_1 \in
  \operatorname{Spec}(\widehat{R})$ existiert mit $P_1 \cap R =
  \mathfrak{p}_1$. Nach (\ref{4.4}) muss $P_1 \cap R$ regulär sein.\\
  Sei im \emph{2. Schritt} $N$ wie in der Voraussetzung und $N_1 \subset E$
  wie in (\ref{4.1},~b) der einzige einfach-teilbare Untermodul von $E$ mit
  $N \sim N_1$. Dann erfüllt auch $N_1$ die Bedingung (*), und nach dem
  1.~Schritt ist $\Ann_R(N) = \Ann_R(N_1)$ ein maximales Element in
  $\Ass(R)$.
\end{Beweis}

\section*{Appendix}

Sei weiter $R$ noethersch und lokal. Vertauscht man die Begriffe teilbar
bzw. torsionsfrei aus den bisherigen Überlegungen mit radikalvoll bzw.
sockelfrei, so erhält man strengere Aussagen, die einige Ergebnisse aus
\cite{005} vervollständigen.

Die erste Beobachtung ist, dass auch jeder einfach-sockelfreie $R$-Modul $M$
prim ist (denn für jedes $r \in R$ ist $M/M[r] \cong rM$ sockelfrei, also
$M[r] = 0$ oder $M[r] = M$), entsprechend jeder einfach-radikalvolle
$R$-Modul $N$ koprim.\\
Die zweite Beobachtung ist, dass die Aussage "`$\Ann_R(M)$ ist ein maximales
Element in $\Ass(R)$"' durch "`$\Ann_R(M)$ ist ein Primideal der Kohöhe~1"'
zu ersetzen ist, genauer:

\bigskip
\noindent
\textbf{Lemma} {\itshape
  Für ein Primideal $\mathfrak{p} \in \operatorname{Spec}(R)$ sind
  äquivalent:
  \begin{iii}
    \item $\kappa(\mathfrak{p})$ ist als $R$-Modul einfach-sockelfrei.
    \item $\kappa(\mathfrak{p})$ ist als $R$-Modul einfach-radikalvoll.
    \item $\dim(R/\mathfrak{p}) = 1$.
  \end{iii}
}

\begin{Beweis}
  Die ersten beiden Aussagen sind äquivalent damit, dass
  $\kappa(\mathfrak{p})$ als $R/\mathfrak{p}$-Modul einfach-sockelfrei bzw.
  einfach-radikalvoll ist, so dass wir ab jetzt $\mathfrak{p} = 0$ annehmen
  können und dann $\kappa(\mathfrak{p})$ der Quotientenkörper $K$ von $R$
  ist.\\
  Aus (iii), d.\,h. $\dim(R) = 1$, folgt sockelfrei = torsionsfrei und
  radikalvoll = teilbar, so dass (i) und (ii) klar sind. Bei (i $\to$ iii)
  gilt für jedes Ideal $\mathfrak{a} \neq 0$, dass $K/\mathfrak{a}$
  halbartinsch ist (vgl. \ref{2.1},~a), also $\dim(R) = 1$, und bei (ii
  $\to$ iii) gilt für jedes Primideal $\mathfrak{q} \subsetneqq
  \mathfrak{m}$, dass $R_{\mathfrak{q}}$ radikalvoll, also $R_{\mathfrak{q}}
  = K$, d.\,h. $\mathfrak{q} = 0$ ist, also ebenfalls $\dim(R)=1$.
\end{Beweis}

\medskip
\noindent
\textbf{Satz 1} {\itshape
  Für einen $R$-Modul $M \neq 0$ sind äquivalent:
  \begin{iii}
    \item $M$ ist einfach-sockelfrei.
    \item $M \hookrightarrow \kappa(\mathfrak{p})$ für ein Primideal
      $\mathfrak{p}$ mit $\dim(R/\mathfrak{p}) = 1$.
  \end{iii}
}

\begin{Beweis}
  (i $\to$ ii) Wie in (\ref{2.1},~a) zeigt man für jeden Untermodul $0 \neq
  U \subset M$, dass $M/U$ halbartinsch und $U$ wieder einfach-sockelfrei
  ist. Insbesondere ist $M$ uniform, mit $\Ass(M) = \{\mathfrak{p}\}$, also
  $M \subset E(R/\mathfrak{p})$, und weil auch $R/\mathfrak{p}$
  einfach-sockelfrei ist, $\dim(R/\mathfrak{p}) = 1$. Nach der ersten
  Beobachtung ist $M$ auch prim, also nach (\ref{1.1},~a) $\Ann_R(M) =
  \mathfrak{p}$, $M \subset E(R/\mathfrak{p})[\mathfrak{p}] \cong
  \kappa(\mathfrak{p})$.\\
  (ii $\to$ i) Nach dem Lemma ist $\kappa(\mathfrak{p})$ einfach-sockelfrei,
  also auch der Untermodul $M$.
\end{Beweis}

\medskip
\noindent
\textbf{Bemerkung} Satz~1 und (\ref{2.2}) lassen sich für jede
Gabriel-Topologie $\mathfrak{G}$ auf $R$ formulieren: Ein $R$-Modul $M \neq
0$ heiße \emph{einfach-$\mathfrak{G}$-torsionsfrei}, wenn $M$
$\mathfrak{G}$-torsionsfrei ist (d.\,h. $\Ass(M) \cap \mathfrak{G} =
\varnothing$) und wenn aus $0 \neq U \subsetneqq M$ stets folgt, dass $M/U$
nicht $\mathfrak{G}$-torsionsfrei ist. Ist nun $\mathfrak{G}^*$ die Menge
aller maximalen Elemente in $\operatorname{Spec}(R)\setminus \mathfrak{G}$,
kann man wie in (\ref{2.1}) und (\ref{2.2}) zeigen, dass $M \neq 0$ genau
dann einfach-$\mathfrak{G}$-torsionsfrei ist, wenn es eine Einbettung $M
\hookrightarrow \kappa(\mathfrak{p})$ gibt mit $\mathfrak{p} \in
\mathfrak{G}^*$.\\
$\mathfrak{G}_1 = \{\mathfrak{a} \subset R \,|\, \mathfrak{a}\text{ ist
  regulär}\}$ liefert dann (\ref{2.2}) (denn $\mathfrak{G}_1^*$ ist die Menge
aller maximalen Elemente in $\Ass(R)$), und $\mathfrak{G}_2 = \{\mathfrak{a}
\subset R \,|\, R/\mathfrak{a}\text{ ist artinsch}\}$ liefert Satz~1 (denn
$\mathfrak{G}_2^*$ ist die Menge aller Primideale der Kohöhe~1).\hfill$\Box$

\bigskip
Jeder einfach-radikalvolle $R$-Modul ist nach (\cite{005} Lemma~4.1)
entweder sockelfrei (und dann uniform) oder artinsch (und dann kouniform):

\bigskip
\noindent
\textbf{Satz 2} {\itshape
  Für einen $R$-Modul $N$ sind äquivalent:
  \begin{iii}
    \item $N$ ist einfach-radikalvoll und sockelfrei.
    \item $N$ ist radikalvoll und einfach-sockelfrei.
    \item $N \cong \kappa(\mathfrak{p})$ für ein Primideal $\mathfrak{p}$
      mit $\dim(R/\mathfrak{p}) = 1$.
  \end{iii}
}

\begin{Beweis}
  (iii $\to$ i) und (iii $\to$ ii) folgen sofort aus dem Lemma, ebenso (ii
  $\to$ iii) mit Satz~1. Bleibt (i $\to$ iii) zu zeigen: Nach (\cite{005}
  Lemma~4.1) ist $\Ann_R(N)$ ein Primideal, sagen wir $\mathfrak{p}$, mit
  $\dim(R/\mathfrak{p}) = 1$. Es folgt $\Ass(N) = \{\mathfrak{p}\}$, weil
    $N$ uniform ist genauer $N \subset E(R/\mathfrak{p})$, also in $N
    \subset \kappa(\mathfrak{p})$ nach dem Lemma Gleichheit. 
\end{Beweis}

\noindent
\textbf{Satz 3} {\itshape
  Für einen $R$-Modul $N$ sind äquivalent:
  \begin{iii}
  \item $N$ ist einfach-radikalvoll und artinsch.
  \item $N \sim \Ann_E(P)$ für ein $P \in \operatorname{Spec}(\widehat{R})$
    mit $\dim(\widehat{R}/P) = 1$.
  \end{iii}
}

\begin{Beweis}
  (i $\to$ ii) Jeder Untermodul $V \subsetneqq N$ ist endlich erzeugt, so
  dass auch $N/V$ einfach-radikalvoll ist. Und weil $N$ kouniform ist, folgt
  aus $\mathfrak{m}^nV = 0$ und $\mathfrak{m}^nN = N$, dass $rN=N$ ist für
  ein $r \in \mathfrak{m}^n$, also $rV = 0$, $N/V \twoheadrightarrow N/N[r]
  \cong N$, $N \sim N/V$. Für irgendein $0 \neq f \colon N \to E$ definieren
  wir jetzt $V = \operatorname{Ke} f \subset N$ und $N_1 = \operatorname{Bi}
  f \subset E$. Dann ist $N \sim N_1$ und $\Ann_{\widehat{R}}(N_1)$ nach
  (\cite{005} Lemma~4.1) ein Primideal in $\widehat{R}$, sagen wir $P$, mit
  $\dim(\widehat{R}/P) = 1$. Wegen $\Ann_E(P) = N_1$ folgt alles.\\
  (ii $\to$ i) Der $\widehat{R}$-Modul $\widehat{R}/P$ ist
  einfach-sockelfrei, also via Matlis-Dualität der $\widehat{R}$-Modul
  $\Hom_{\widehat{R}}(\widehat{R}/P,E) \cong \Ann_E(P)$ einfach-radikalvoll.
  Natürlich ist dann $\Ann_E(P)$ auch als $R$-Modul einfach-radikalvoll (und
  artinsch), also auch der Faktormodul $N$.
\end{Beweis}

\end{document}